\documentclass[num-refs]{wiley-article}

\usepackage{siunitx}
 \usepackage[utf8]{inputenc}
\usepackage{hyperref}
\hypersetup{
    colorlinks=true,
    linkcolor=black,
    filecolor=magenta, urlcolor=cyan,
    citecolor=blue,
    } 

\usepackage{graphicx,amsmath,url} 
\usepackage{float}
\usepackage{indentfirst}

\usepackage{comment}

 \usepackage{amssymb} 

\newcommand\abs[1]{\ensuremath{\lvert#1\rvert}}

\newcommand{\sgn}{\operatorname{sgn}}

\papertype{Original Article}
\paperfield{Journal Section}

\title{Global Output-Feedback Extremum Seeking Control with Source Seeking Experiments}

\author[1\authfn{1}]{Nerito Oliveira Aminde}
\author[2\authfn{2}]{Tiago Roux Oliveira}
\author[3\authfn{3}]{Liu~Hsu}

\affil[1]{Federal Center for Technological Education "Celso Suckow da Fonseca"
 (CEFET/RJ), Angra dos Reis, Rio de Janeiro -- RJ, 23953-030, Brazil}
 \affil[2]{State University of Rio de Janeiro (UERJ), Rio de Janeiro -- RJ, 20550-900, Brazil}
 \affil[3]{Federal University of Rio de Janeiro (COPPE/UFRJ), Rio de Janeiro -- RJ, 21941-972, Brazil}

\corraddress{Nerito Oliveira Aminde}
\corremail{nerito.aminde@cefet-rj.br}

\presentadd[\authfn{1}]{Department of Electrical Engineering, CEFET/RJ}
\presentadd[\authfn{2}]{Department of Electronics and Telecommunication Engineering, UERJ}
 \presentadd[\authfn{3}]{ Program of Electrical Engineering, COPPE/UFRJ}

\fundinginfo{CAPES, CNPq, FAPERJ, ISDB-Maputo and CINOP Global}

\runningauthor{Aminde et al.}

\begin{document}

\maketitle

\begin{abstract}
This paper discusses the design of an extremum seeking controller that relies on a monitoring function for a class of SISO uncertain nonlinear systems characterized by arbitrary and uncertain relative degree. Our demonstration illustrates the feasibility of achieving an arbitrarily small proximity to the desired optimal point through output feedback. The core concept involves integrating a monitoring function with a norm state observer for the unitary relative degree case and its expansion to arbitrary relative degrees by means of the employment of a time-scaling technique. Significantly, our proposed scheme attains the extremum of an unknown nonlinear mapping across the entire domain of initial conditions, ensuring global convergence and stability for the real-time optimization algorithm. Furthermore, we provide tuning rules to ensure convergence to the global maximum in the presence of local extrema. To validate the effectiveness of the proposed approach, we present a numerical example and apply it to a source-seeking problem involving a cart-track linear positioning servomechanism. Notably, the cart lacks the ability to sense its velocity or the source's position, but can detect the source of a light signal of unknown concentration field.

\begin{small}
\end{small}
\end{abstract}

\section{Introduction}

Extremum seeking control (ESC) is a control system used to determine in real-time the extremum (maximum or minimum) of an unknown nonlinear mapping  \citep{NMMM:10}.  
ESC has become a key area in control theory due to increasing need to optimize plant operation in order to reduce operating costs and meet product specifications \citep{GAD:04,CWPF:20}.
In \cite{NMMM:10} and the cited references, numerous applications are discussed, including but not limited to the design of antilock braking systems, autonomous vehicles, mobile robots, internal combustion engines, process control, particle accelerators, and source seeking. Recently, an extension of the scope of applicability of the ESC through delays and PDEs was proposed in \citep{OK:22}, with practical engineering scenarios including neuromuscular electrical stimulation \citep{POPF2020}, non-cooperative games \citep{ORKBEK:19}, biological reactors \citep{RPFT:19,OFKK2020}, oil-drilling systems \citep{OK:20}, and flow-traffic control for urban mobility \citep{YSOK:2020}.  

The prevalent algorithms for unconstrained optimization typically rely on the derivative or gradient of the objective function. Nonetheless, in numerous control problem applications, including those outlined earlier, crucial components such as the plant model, gradient information, and the cost function for optimization are not readily accessible online. Furthermore, "gradient sensors" have a tendency to magnify noise and encounter instability issues, particularly at higher frequencies \cite{OALE:07}.  

Classical ESC methods employ a high-pass filter in the system output along with a small sinusoidal perturbation (dither signal) technique to estimate the gradient of the cost function. This approach is known for its simplicity and rapid adaptation \cite{AK:03,KW:00,LGK:18}. Nevertheless, only local stability properties could be assured when assuming full-state measurement.
In \cite{NMMM:10,TNM:06,TNMA:09}, 
under the identical assumption, semi-global practical convergence was achieved, but with a decrease in the convergence rate within the domain of attraction.

In \cite{POL:12,KU:74,POA:03}, the extremum seeking could be seen as a nonlinear control challenge featuring a state-dependent high-frequency gain (HFG), which alters direction (referred to as control direction) around the optimal point of interest. In addressing this, a method employing sliding mode control for tracking uncertain plants with an unknown control direction was introduced in 
\cite{POLN:07,HYCL:08, HOP:10} using an algorithm of switching based on a monitoring function to the output error. In  \cite{HOP:11}, it was conjectured that the lack of
robustness of such monitoring scheme  with respect to recurrent changes of HFG sign would preclude the monitoring function  to be directly applied to extremum seeking control. 

In this paper, a novel monitoring function is proposed in order to show that the output-feedback tracking controller proposed in \cite{POLN:07,HYCL:08,HOP:10}, indeed can also be applied to ESC of a class of uncertain nonlinear systems (of unitary and arbitrary relative degrees), while global convergence properties of the search algorithm are also guaranteed without affecting its rate of convergence. Moreover, the proposed algorithm can be reconfigured to achieve a global extremum point in the presence of local extrema.

On the other hand, the control of systems with uncertain relative degrees is a
challenging task. There are few results in the adaptive and sliding mode control literature
\cite{M:1985,M:1987,M:2000,BPU:2008} which consider the model reference and stabilization problems, but not for extremum
seeking. Even in those cases the processes are restricted to an
uncertainty not exceeding 2, \emph{i.e.}, the uncertain relative
degree may be $n^*$, $n^*+1$ or $n^*+2$.
As a further contribution of this paper, a generalization is achieved to include more general
dynamics with arbitrary and uncertain relative degrees, without
using differentiators to compensate them \cite{Lev:03,AFL:2012,FSEY:2008,BFU:1998}.
The theoretical contribution of this paper is to develop a
time-scaling procedure in order to reduce the order of the system
dynamics, and consequently, to allow the analysis and control
design for enlarged uncertainty in the relative degree. 

In particular, research in applications that use autonomous vehicles are wide,
varied, and constantly growing. Notably, there is a burgeoning interest in researching vehicles that operate without access to precise position information. These vehicles must navigate and perform
a desired task without the use of GPS (Global Positioning System)
or inertial navigation \cite{G:2011,ZAGSK:2007,CGSK:2009,ZDZ:2023}. A recent solution to this problem was given in the context of
source seeking
\cite{ZAGSK:2007,CGSK:2009,GK:2011,FO:2009,CKKXK:2009}. This
problem has been intensively explored in the robotics literature
where results of extremum seeking control
\cite{KW:00,AK:03,POA:03,TNM:06,TNMA:09,NMMM:10} are applied to
autonomous mobile robots with the objective of localizing the
source of an unknown, nonlinear, signal field. For environments where position information is unavailable, the extremum seeking method is applied to autonomous vehicles as a means of navigating to find the source of some signal which the vehicles can measure locally. The signal is at maximum intensity at the source and decreases with distance away from the source. The source signal to be detected may be electromagnetic, acoustic, a temperature intensity or the concentration of a biological/chemical
agent \cite{GK:2011,FO:2009,CKKXK:2009}.

Finally, we conduct experiments to validate 
the proposed theoretical 
results. 
We employ the new ESC scheme based on monitoring functions in a
light-source seeking scenario for a one-dimension optimization problem,
\emph{i.e.}, both input and output signals are scalar functions.
The theoretical results for basic extremum seeking are applied to
control a linear servo-positioning system to perform localization
and tracking of a
light source, without its position/velocity information. Preliminary conference versions were presented in \cite{AOH:13,AOH:14}.

\section{Preliminaries and problem formulation} \label{section2}

Throughout the paper, the Euclidean norm of a vector $x$ and the corresponding induced norm of a matrix $A$ are denoted by $\|x\|$ and $\|A\|$, respectively. The term $\pi_{i}(t)$ stands for any exponential decaying function, such that $\abs{ \pi_{i} (t) } \leq R e^{-\beta t},\forall t$, and some positive scalars $R$ and $\beta$. Class  $\mathcal{K}$ and $\mathcal{K}_\infty$  functions are defined as in \cite{K:02b}.
From a technical standpoint, the theoretical results obtained in this paper are based on Filippov's definition for solution of differential equations with discontinuous right-hand sides \cite{F:64}. 

Consider  the following nonlinear uncertain system:
\begin{eqnarray}             
	\dot{x}&=&f(x,t)+g(x,t)u                                                                               \label{eq1}  \\                                                                       
	z&=&h(x,t)                                                                                                 \label{outupunmeasured}
\end{eqnarray}
in cascade with a static subsystem
\begin{align}
y=\Phi(z)\,,                                                                                                \label{outputmensured}
\end{align}
where $u \in \mathbb{R}$ is the control input,  $x \in \mathbb{R}^n$ is the state vector, $z \in \mathbb{R}$  and  $y \in \mathbb{R}$ are measured outputs of the first
subsystem and of the static subsystem, respectively. In order to assure existence and forward uniqueness
of solutions, the nonlinear uncertain functions $f,g$ and $h$ are locally Lipschitz continuous in $x$, piecewise continuous in  $t$
and sufficiently smooth. Without loss of generality, we assume that the
initial time is $t=0$. For each solution of (\ref{eq1})
there exists a maximal time interval of definition given by
$[0,t_M)$, where $t_M$ may be finite or infinite.
The control objective of ESC is not ``stabilization'' or ``tracking'', but is ``real-time optimization'' \cite{TNMA:09}.  However, the ESC problem can be reformulated as a tracking problem in which the control direction is unknown \cite{POL:12}. We wish to find an output-feedback control law $u$ so that, from any initial condition, the system  is steered to reach the extremum point $y^*$ and remain on such point thereafter, as close as possible. Without lost of generality, we only address the maximum seeking problem.

The system (\ref{eq1})--(\ref{outputmensured})  can be rewritten in the normal form as follows: 
\begin{eqnarray}
  \dot{\eta} &=& \phi_0(\eta,z,t)\,,                                              \label{plant_inverse2}\\
  \dot{z} &=& \phi_1(\eta,z,t)+ \phi_2(\eta,z,t) u\,,                     \label{plant_extern2}\\
  y&=&\Phi(z)\,,                                                                          \label{saidameasured}
  \end{eqnarray}
 with state $x:=\left[\eta^T \ z\right]^T$, $\eta \in \mathbb{R}^{n-1}$ and $z \in \mathbb{R}$, and uncertain nonlinear functions $\phi_0: \mathbb{R}^{n-1} \!\times\!
\mathbb{R}\!\times\! \mathbb{R}^{+}\to\! \mathbb{R}^{n-1}$ and $\phi_1,\phi_2: \mathbb{R}^{n-1} \!\times\! \mathbb{R}\!\times\! \mathbb{R}^{+} \!\to\! \mathbb{R}$.

The state $\eta$ of the $\eta$-subsystem, referred to as an ``inverse system'', is not available for feedback.

 With respect to the controlled plant, we assume the following assumptions:

$\bf{(A1)}$ (\emph{On the uncertainties}): All the uncertain plant parameters belong to a compact set $\Omega$.

This assumption is necessary to obtain the uncertainty bounds for control design. 

$\bf{(A2)}$ (\emph{Relative degree one}): The uncertain function $\phi_2(\eta,z,t)$ is bounded away from  
zero, \emph{i.e.}, $$0 < \underline{\phi}_{2} \leq |\phi_2|\,, \quad \forall t \in [0,t_M)\,,$$
where the constant lower bound $\underline{\phi}_{2}$ is known.

According to $\bf{(A2)}$,  the subsystem
(\ref{plant_inverse2})--(\ref{plant_extern2}) has relative degree one \emph{w.r.t.} $z$ since $\phi_2 \neq 0$.
It restricts us to the case of relative degree one, which is the simplest case amenable by pure Lyapunov design.

By using the notation  $\Phi'(z):=\frac{d \Phi(z)}{dz}$ and $\Phi''(z):=\frac{d^2 \Phi(z)}{dz^2}$, we consider that (in $\Omega$):

$\bf{(A3)}$ (\emph{Cost Function}):   There exists a unique point $z^*\in\mathbb{R}$ such that $y^*=\Phi(z^*)$ is the
extremum (maximum) of $\Phi(z)$: $\mathbb{R} \rightarrow \mathbb{R}$, satisfying
\begin{align}
\Phi'(z^*)=0,\,\,\,\Phi''(z^*)<0  \nonumber  \\                                                    
\Phi(z^*)>\Phi(z), \,\, \forall z\in\mathbb{R},\,\,z \neq z^*                                      \nonumber
\end{align}
and for any given $\Delta>0$, there exists a constant $L_\Phi(\Delta)>0$ such that 
\begin{align}
L_\Phi(\Delta)\leq|\Phi'(z)|,\,\forall z\notin \mathcal{D}_\Delta:=\{z:|z-z^*|<\Delta/2\},  \notag 
\end{align}
where $\mathcal{D}_\Delta$ is called $\Delta$-vicinity of  $z^*$ and $\Delta$ can be made arbitrary small by allowing a smaller $L_{\Phi}$.

From (\ref{plant_extern2}) and (\ref{saidameasured}), the first time derivative of the output $y$ is given by
\begin{align}
\dot{y} = \Phi' \phi_1 + k_p u\,,\label{eq:ydynamics}
\end{align}
where the plant high frequency gain (HFG), denoted by $k_p$, is the coefficient of $u$ which appears in the first time derivative of the output $y$ and it is given by
\begin{align}
k_p=\Phi'  \phi_2\,.                                                                                              \label{kpneras}
\end{align}

As in \cite{POL:12}, the $\sgn(k_p)$ is also called \textit{control direction}. The assumption $(\bf{A3})$  leads us to consider a nonlinear control system with a state dependent HFG which changes sign around the optimum point of interest in a continuous way. 

From (\ref{kpneras}), $(\bf{A2})$ and $(\bf{A3})$, $k_p$ satisfies $(\forall z \notin \mathcal{D}_{\Delta})$
\begin{align}
0<\underline{k}_{p}\leq |k_p|                                                                               \label{kpbarneras}
\end{align}
where the lower bound $\underline{k}_{p} \leq \underline{\phi}_{2} L_\Phi$ is a constant.

$\bf{(A4)}$ (\emph{Norm observability}): The inverse system (\ref{plant_inverse2}) admits a known first order norm observer  of the form:
\begin{equation}
\dot{\bar{\eta}} = - \lambda_0 \bar{\eta} + \varphi_0(z,t)\,,                              \label{normobsgeneric}
\end{equation}
with $z$ in (\ref{plant_extern2}), input $\varphi_0(z,t)$ and output $\bar{\eta}$ such that: (i) $\lambda_0\!>\!0$ is a
constant; (ii) $\varphi_0(z,t)$ is a non-negative function continuous in $z$, piecewise continuous and upper-bounded in $t$;
and (iii) for each initial states $\eta(0)$ and $\bar{\eta}(0)$
\begin{equation}                                                             
\|\eta\| \leq |\bar{\eta}| + \pi_0\,, \quad \forall t\in[0,t_M)\,,                                    \label{etabound}
\end{equation}
where
$\pi_0\!:=\!\Psi_0(|\bar{\eta}(0)|\!+\!\|\eta(0)\|)e^{-\lambda_0 t}\!$ and $\Psi_0\in\mathcal{K}$. 

It is well known that, if the inverse system (\ref{plant_inverse2}) is input-to-state-stable (ISS) w.r.t. $z$, then it admits such norm observer and the plant is minimum-phase \cite{KSW:01}. More examples of nonlinear systems which satisfies such assumption are given in \cite{POL:12}.

In order to obtain a norm bound for the term $ \Phi'\phi_1$ in (\ref{eq:ydynamics}), we additionally assume that:

$\bf {(A5)}$ (\emph{Domination Functions}): There exist {\em known} functions $\bar{\Phi},\alpha_1\in\mathcal{K}_\infty$, with $\alpha_1$ locally Lipschitz, a {\em known} non-negative function $\varphi_1(z,t)$ continuous in $z$, piecewise continuous and upper bounded in $t$ such that \\$|\phi_1(\eta,z,t)|\!\leq\!\alpha_1(\|\eta\|)\!+\!\varphi_1(z,t)$ and $|\Phi'|\!\leq\!\bar{\Phi}(|z|)$.

Note that the Assumption $(\bf{A5})$ is not restrictive, since $\Phi'$ is assumed to be smooth and $\phi_1$ is locally Lipschitz continuous in $\eta$ and in $z$. Furthermore, the domination functions $\alpha_1$ and $\varphi_1$ impose stringent growth condition only \emph{w.r.t.} the time-dependence. Thus, polynomial nonlinearities in $\eta$ and $z$ can be coped with.

\section[Output-feedback extremum seeking controller via monitoring function]{Output-feedback extremum seeking controller via \\monitoring function} \label{section3}

\begin{figure}[!htb]
\begin{center}
\includegraphics[width=.55\textwidth]{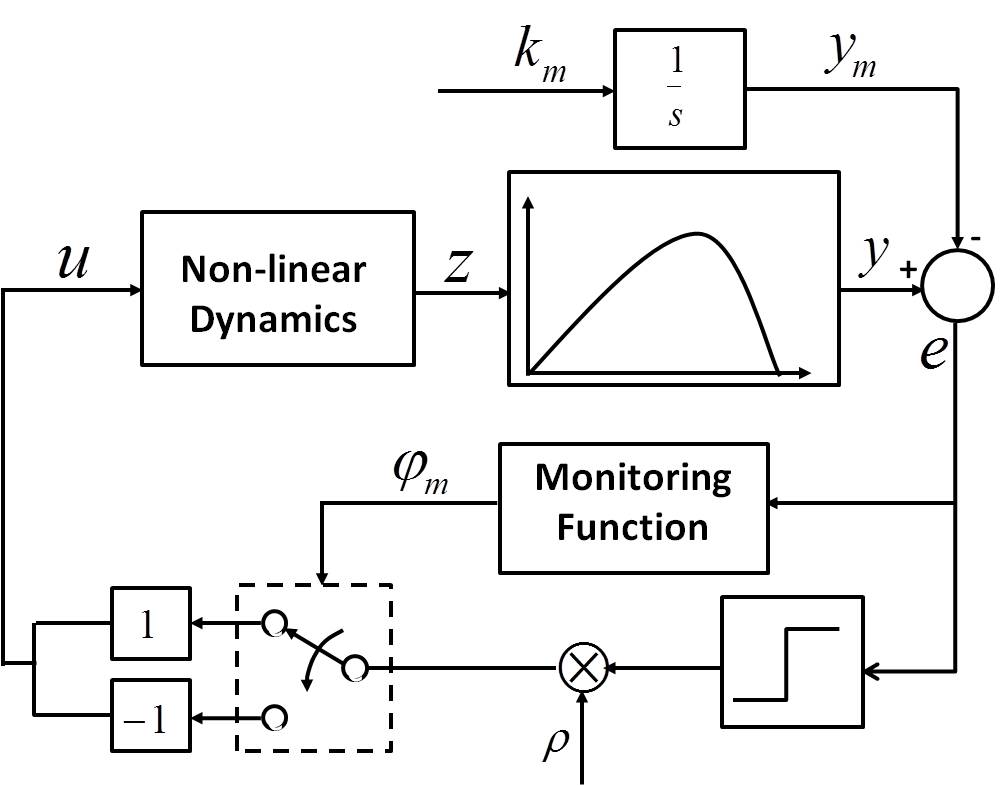}
\caption{Extremum seeking controller using a monitoring function.}
\label{diagram_dynamic}
\end{center}
\end{figure}

The proposed  output-feedback  ESC based on monitoring function is represented in Fig. \ref{diagram_dynamic}. The control law to plants with unknown HFG is defined as in \cite{POLN:07,HYCL:08}
\begin{equation}
 u=\left\{
	\begin{array}{rcl}
	u^+&=&-\rho\sgn(e)\,, \quad t\in T^+\,,\\
	u^-&=& \rho\sgn(e)\,, \quad t\in T^-\,,
	\end{array}\right.        
	\label{control_law}
 \end{equation}
where the monitoring function is used to decide when $u$ should be switched from $u^+$ to $u^-$ and vice versa. In (\ref{control_law}), $\rho$ is the modulation function to be designed later on and the sets $T^+$ and $T^-$ satisfy $T^+~\cap~T^- = \emptyset$ and $T^+~\cup~T^- = [0,t_M)$.

The  tracking error $e$ is given by the following
\begin{equation}
e(t) = y(t) - y_m(t)\,,                                                                                        \label{erro}
\end{equation}
where  $y_m\in\mathbb{R}$ is a simple ramp time function generated by the reference model 
\begin{equation}
\dot{y}_m = k_m,  \quad y_m(0)=y_{m0} \,,                                                 \label{modref} 
\end{equation}
where $k_m$ and $y_{m0}$ are design constants. In order to avoid an unbounded reference signal $y_m(t)$ in the controller, one can saturate the model output at a rough known norm upper bound of $y^*$ without affecting the control performance.

The modulation function $\rho$ will be designed so that $y(t)$ tracks $y_m(t)$ as long as possible.
In this way, $y$ is forced to achieve the vicinity of the maximum $y^*=\Phi(z^*)$ and remains close to the optimum value $y^*$.
To this end, we have to propose $\rho$ such that the output tracking error $e$ tends to zero in finite time at least outside the $\Delta$-vicinity, that is, in the neighborhood of the maximizer $z^*$.

Thus, it is straightforward to 
conclude that $y=\Phi(z)$ tries to track $y_m$ (and consequently $y$ must approach the maximum at $y^*$) as long as $y$ remains away from a small vicinity of $y^*$ where the HFG is away from zero. In contrast, once $y$ reaches the vicinity of $y^*$, the HFG will approach zero and thus controllability is lost. Consequently, tracking of $y_m$ will cease. But then the neighborhood of the optimum point is already achieved as desired. Our control strategy will guarantee that $y$ will remain close to $y^*$ thereafter. It is apparent that the convergence rate of $z$ to the $\Delta$-vicinity $\mathcal{D}_\Delta$ defined in $\bf{(A3)}$ is a function of $\rho$. 
Although $\mathcal{D}_\Delta$ is not positively invariant, after reaching $\mathcal{D}_\Delta$, it will be shown that $z$ will remain close to $z^*$ where the maximum takes place. It does not imply that
$z(t)$ remains in $\mathcal{D}_\Delta$, $\forall t$. However, as shown later on in Theorem~\ref{theorem:01}, one can guarantee that $y$ remains close to the optimum value $y^*$.

\subsection{Error Dynamics}

From (\ref{eq:ydynamics}), (\ref{erro}) and (\ref{modref}), by adding and subtracting $\lambda e$ the time derivative of the error $e$ one has 
\begin{align}
\dot{e}&=\Phi' \phi_1+k_pu-k_m+\lambda e-\lambda e\,, \\
\dot{e}&=-\lambda e+k_p(u+d_e)\,,                                                  \label{edoerro}
\end{align}
where $\lambda>0$ is an appropriate constant and

\begin{align}
d_e:=\frac{1}{k_p}[{\Phi' \phi_1} - {k_m} + {\lambda e}] \,.         \label{moduerro}
\end{align}

In \cite{HOP:10}, it is shown that  if the control law 
\begin{align}
	u=-\sgn(k_p)\rho\sgn(e)
\end{align}

 was used with modulation function
$\rho$ satisfying
\begin{align}
\rho\geq|d_e|+\delta,                                                                     \label{funcmodgen}
\end{align}
{\em modulo} exponential decaying terms, then by using the comparison lemma \cite{F:64}, one has $\forall t \in
[t_i,t_M)$:
 \begin{equation}
|e(t)|\leq  \zeta(t)\,, \quad \zeta(t):=|e(t_i)| e^{-\lambda
(t- t_i)}\!+\!\pi_1\,,\label{eq:boundMIMOKpconhecido}
\end{equation}
where $\pi_1\!:=\!\Psi_1(|\bar{\eta}(0)|+\|\eta(0)\|)e^{-\lambda_1
t}$, $0\!<\!\lambda_1\!<\!\min\{\lambda_0,\lambda\}$ and
$\Psi_1\!\in\!\mathcal{K}$.

On the other hand, if inequality (\ref{funcmodgen}) were verified taking into account the exponential decaying terms, the upper bound  (\ref{eq:boundMIMOKpconhecido}) would be modified to
 \begin{equation}
|e(t)|\leq  \zeta(t)\,, \quad \zeta(t):=|e(t_i)| e^{-\lambda
(t- t_i)}\!\,.\label{eq:boundMIMOKpconhecido2}
\end{equation}

The major problem is that $\sgn(k_p)$ is unknown in both cases, thus we cannot
implement it.
In what follows, a switching scheme based on monitoring
function is developed to cope with the lack of knowledge of the
control direction outside the $\Delta$-vicinity.

\subsection{Monitoring Function Design}
\label{monitsection}
The detailed description of monitoring function can be found in \cite{HOP:10}. Here, only a brief description of how it works is given.
Reminding that inequality (\ref{eq:boundMIMOKpconhecido2}) holds when the control direction is correct, it seems natural to use
$\zeta$ in (\ref{eq:boundMIMOKpconhecido2}) as a benchmark to decide whether a switching of $u$ in (\ref{control_law})  from $u^+$ to $u^-$ (or $u^-$ to $u^+$) is  needed, \emph{i.e.}, the switching occurs only when (\ref{eq:boundMIMOKpconhecido2}) is violated.

Therefore, consider the following function,
\begin{align} 
\varphi_k(t)= |e(t_k)|e^{-\lambda(t-t_k)}+ r,                              \label{eqmonitpart}
\end{align}   
where $t_k$ is the switching time and 
$r$ is any arbitrary small constant. The monitoring function $\varphi_m$ can be defined as
\begin{align}
\varphi_m(t):=\varphi_k(t), \quad \forall t \in [t_k,t_{k+1}) \subset [0,t_M).                 \label{eqmonit}
\end{align}
Note that, from (\ref{eqmonitpart}) and (\ref{eqmonit}), one has $\abs{e(t)} < \abs{\varphi_k(t)}$ at $t=t_k$. Hence, $t_k$ is defined as the time instant when the monitoring function $\varphi_m(t)$ meets $\abs{e(t)}$, that is,
\begin{equation}
\label{eq:tk_nstar>1}
 t_{k+1}:=
    \begin{cases}
     \min\{t > t_{k}: |e(t)| = \varphi_k(t)\}, & \text{if it exists}\,, \\
     t_M, & \text{otherwise}\,,
    \end{cases}
\end{equation}
where $k \in \{0,1,\ldots\}$ and $t_0\!:=\!0$ (see
Fig.~\ref{funmonit_error}).

\begin{figure}[!htb]
\begin{center}
\includegraphics[width=.6\textwidth]{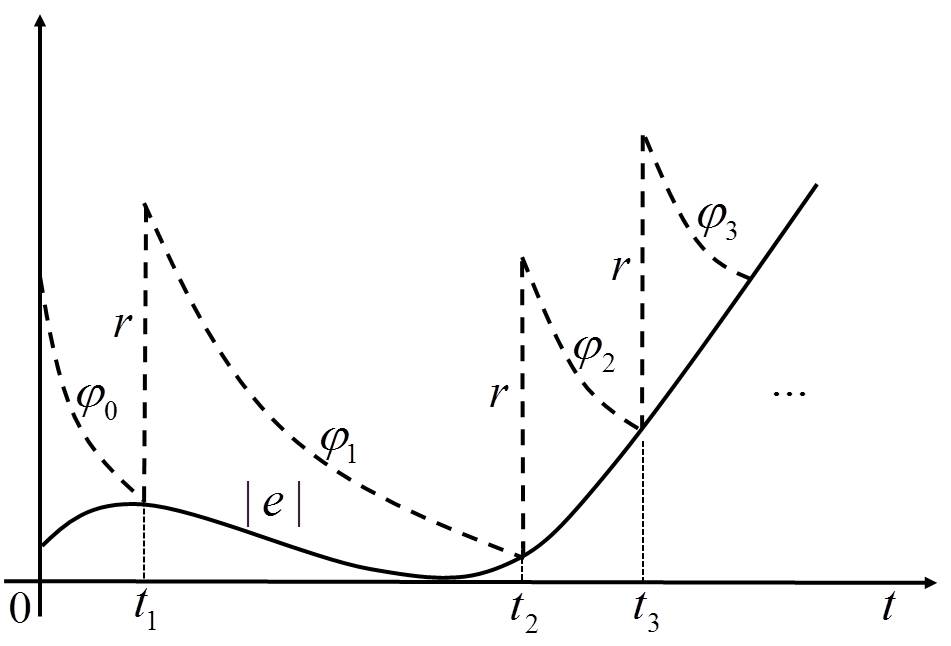}
\caption{The trajectories of $\varphi_m(t)$ and $|e(t)|$.}
\label{fig1monit}
\end{center}
\end{figure}

The following inequality is directly obtained from
(\ref{eqmonit})
\begin{equation}
|e(t)| \leq \varphi_m(t), \quad \forall t \in [0,t_M) \,.
\label{eq:bound_phim}
\end{equation}
Fig.~\ref{fig1monit} illustrates the tracking
error norm $|e|$ as well as the monitoring function $\varphi_m$.

\begin{remark} [Main Changes in the Monitoring Function:]
\emph{The main difference with respect to  \cite{POLN:07,HYCL:08, HOP:10} is that the monitoring function previously proposed  was based on the upper bound (\ref{eq:boundMIMOKpconhecido}) such that $\varphi_k$ in (\ref{eqmonitpart}) is replaced by
\begin{equation}
\varphi_k(t) = |e(t_k)| e^{-\lambda(t- t_k)} + a(k)
e^{- t/a(k)}\,, \label{monitor_mimo}
\end{equation}
where  $a(k)$ is any positive monotonically increasing
unbounded sequence. The purpose of the second term in (\ref{monitor_mimo}) is to dominate the term $\pi_1$ which is not available for measurement. In our ESC application it might be a problem since, as will be shown, the ultimate residual set of the proposed algorithm around the maximum $y^*$ is dependent on the values to which the monitoring function converges. Since (\ref{monitor_mimo}) can assume arbitrarily large values which eventually result in large transients as a result of a change of the control direction and repetitive switching ($k\to+\infty$). In contrast, this does not occur with the new definition given in (\ref{eqmonitpart}) whereby the ultimate residual set can be fixed in the order $\mathcal{O}(r)$}.
\end{remark}

\subsection{Modulation Function Design}

The following auxiliary available upper bounds provide one possible design for the modulation function so that (\ref{funcmodgen}) holds, and they are obtained, taking into account exponentially decaying terms, by using the norm observer given in \textbf{(A4)} and the bounding functions given in \textbf{(A5)}.

From \textbf{(A5)} and (\ref{etabound}), one has $|\phi_1| \leq \alpha_1(|\bar{\eta}|+\pi_0) + \varphi_1(z,t)$. Now, note that $\psi(a\!+\!b)\!\leq\!\psi(2a)\!+\!\psi(2b)$, $\forall a,b\!\geq\!0$ and $\forall \psi\!\in\!\mathcal{K}_\infty$. Thus, since $\alpha_1\!\in\!\mathcal{K}_\infty$ one can write
$\alpha_1(|\bar{\eta}|+\pi_0) \leq \alpha_1(2|\bar{\eta}|)+\alpha_1(2\pi_0)$ and $$|\phi_1| \leq \alpha_1(2|\bar{\eta}|)+\alpha_1(2\pi_0) +\varphi_1(z,t)\,.$$

From (\ref{etabound}), $\pi_0\!:=\!\Psi_0(|\bar{\eta}(0)|\!+\!\|\eta(0)\|)e^{-\lambda_0 t}$ is uniformly bounded. Thus, since $\alpha_1$ is assumed locally Lipschitz in \textbf{(A5)}, one can obtain a valid linear upper bound for $\alpha_1$ such that
$$\alpha_1(2 \pi_0) \leq 2 k_1  \pi_0 = 2 k_1 \Psi_0(|\bar{\eta}(0)|\!+\!\|\eta(0)\|)e^{-\lambda_0 t}\,,$$
where $k_1$ is a positive constant depending on the Lipschitz constant of $\alpha_1$. Then, defining 
\begin{equation}\label{cocozao}
\bar{\phi}_1:=\alpha_1(2|\bar{\eta}|)+\varphi_1(z,t)
\end{equation}
and $\bar\pi_1:=2 k_1 \Psi_0(|\bar{\eta}(0)|\!+\!\|\eta(0)\|)e^{-\lambda_0 t}$, one can also write
\begin{equation}
|\phi_1| \leq \bar{\phi}_1 + \bar\pi_1\,,                                                              \label{boundonphi1}
\end{equation}
where $\bar\pi_1$ decays exponentially.

Furthermore, from (\ref{boundonphi1}), the first term $\phi_1 \Phi'$ of the $y$-dynamics in (\ref{eq:ydynamics}) satisfies $|\phi_1 \Phi'| \leq |\Phi'| |\bar{\phi}_1| + |\Phi'| \bar\pi_1 \leq |\Phi'| |\bar{\phi}_1| + |\Phi^{'}|^2 + \bar\pi_1^2$, where we have used the relationship $|\Phi'| \bar\pi_1 \leq
|\Phi^{'}|^2 + \bar\pi_1^2$. Now, from \textbf{(A5)} the following upper bound holds
\begin{equation}
|\phi_1 \Phi^{'}| \leq \bar{\phi}_1 \bar{\Phi} + \bar{\Phi}^2 + \bar\pi_1^2\,.          \label{defphidag1}
\end{equation}

Remind that, outside the $\Delta$-vicinity, the derivative of the cost function does not vanish $\forall z$. Thus, the lower norm bound
$\underline{k}_p$ for $k_p=\Phi' \phi_2$ given in (\ref{kpbarneras}) holds.

Therefore, one can obtain the following norm bound for $d_e$ defined in (\ref{moduerro}):
\begin{equation}                                                                                          \label{boundond} 
|d_e(t)|\leq \bar{d}_e+\pi_2/\underline{k}_{p}\,, \quad
\bar{d}_e:=(\bar{\phi}_1 \bar{\Phi} + \bar{\Phi}^2+k_m+\lambda |e|)/\underline{k}_{p}\,,
\end{equation}
with the exponential decaying function $\pi_2=\bar\pi_1^2$.

In the proposed scheme,  the following proposition provides one possible modulation function implementation so that (\ref{funcmodgen}) holds. 
\newpage
\begin{proposition}                                                                      \label{prop:01}
Consider the system (\ref{plant_inverse2})--(\ref{saidameasured}), reference model (\ref{modref}) and control law (\ref{control_law}). Outside the $\Delta$-vicinity  $\mathcal{D}_\Delta$, if $\rho$ in (\ref{control_law}) is designed as
\begin{equation}
\rho:=\frac{1}{\underline{k}_p} \left[\bar{\phi}_1 \bar{\Phi} + \bar{\Phi}^2+k_m+\lambda |e|\right] + \Pi(k) + \delta\,,        \label{funcmod}
\end{equation}
 then, while $z \notin \mathcal{D}_\Delta$, one has: \textbf{(a)} the monitoring function switching stops, \textbf{(b)} no finite-time escape occurs in the system signals ($t_M \to +\infty$),  and \textbf{(c)} the error $e$ tends to zero in finite time.  The term $\Pi(k)=a(k)e^{- t/a(k)}$ with $a(k)$ being any positive monotonically increasing
unbounded sequence and $\delta$ is any arbitrary small positive constant. 
\end{proposition}                                                                   
\textbf{Proof:} Outside the $\Delta$-vicinity, suppose by contradiction that $u$ given by (\ref{control_law}) switches without stopping $\forall t\!\in\![0,t_M)$, where $t_M$ may be finite or infinite. Then, $\Pi(k)$ in (\ref{funcmod}) increases unboundedly as $k\!\to\!+\infty$. Thus, there is a finite value $\kappa$ such that for $k\geq\kappa$: (i) the term $\Pi(k)$ will upper bound $\pi_2/\underline{k}_{p}$ in  (\ref{boundond}) and (ii) the control direction is correct. From (i), $\varphi_m(t)\!>\!\zeta(t)$, $\forall t\!\in\![t_{\kappa},t_{\kappa+1})$, with $\zeta$ in (\ref{eq:boundMIMOKpconhecido2}). From (ii), $\zeta$ is a valid upper bound for $|e|$. Hence, no switching will occur after $t=t_{\kappa}$, \emph{i.e.}, $t_{\kappa+1}=t_M$ (see (\ref{eq:tk_nstar>1})), which leads to a contradiction. Therefore, $\varphi_m$ has to stop switching after some finite $k\!=\!N$ and $t_N\!\in\![0,t_M)$, whenever $z \notin \mathcal{D}_\Delta$.

Suppose that we end up with an incorrect control direction estimate. Then, the equation for $e$ can be written as 
$\dot{e}=-\lambda e+ |k_p|[\rho \sgn(e)+d_e]$. In this case, if $\rho$ is defined as (\ref{funcmod}), there exists $t_d < t_M$ such that $e\dot{e}>0, \forall t>t_d$. Hence, $e$ would diverge as $t \to t_M$ for all initial conditions. Therefore, $\sgn(k_p)$ must be correctly estimated at $k=N$.

Now, consider the quadratic function
\begin{align}
	V_e=\frac{e^2}{2}\,.
\end{align}

Then calculating $\dot{V}_e$ along the solutions of (\ref{edoerro}), 
\begin{align}
	\dot{V_e}&=e\dot{e} \\
	\dot{V_e}&=e[-\lambda e+k_p(u+d_e)]\,.
\end{align}
Since the control signal is given by (\ref{control_law}) and the $\sgn(k_p)$ was correctly estimated, the function $\dot{V_e}$ can be rewritten as
\begin{align}
\dot{V_e}=&e[k_p(-\sgn(k_p) \rho \sgn(e) + d_e) - \lambda e]\\
               =&e[-|k_p|\rho \sgn(e)+k_pd_e - \lambda e]\\
							\leq &|e||k_p|[-\rho + |d_e| + |k_p^{-1}\lambda e|].
\end{align}
Since the modulation function $\rho$ satisfies (\ref{funcmod}), the following inequality holds $\dot{V_e} \leq -\delta|k_p||e| <0$, is valid with a constant  $\delta >0$, and the condition $e \dot{e}<0$ (or equivalently $\dot{e}=-\lambda_1\sgn(e)$, for some $\lambda_1>0$), is verified such that $e(t)\to0$ in finite time, while  $z \notin \mathcal{D}_\Delta$. Consequently,  it is not difficult  to conclude that no finite-time escape can occur, \emph{i.e.},
$t_M\!=\!+\infty$. \hfill  $\qed$

\begin{remark}[Modulation Function Reset]
\emph{The term $\Pi(k)$ in (\ref{funcmod}) plays a key role in the domination of the exponential decaying term $\pi_2/\underline{k}_{p}$ in  (\ref{boundond}). It allows that inequality (\ref{funcmodgen}) is satisfied before that  $\pi_2/\underline{k}_{p}$ ultimately vanishes. However, since $\Pi(k)\to+\infty$ as $k\to+\infty$, the modulation function needs a reset mechanism to reinitialize $k$, from time to time, in order to avoid that the controller amplitudes increase to very high values. In particular, if a first order nonlinear system is considered (i.e., the $\eta$-dynamics in  (\ref{plant_inverse2})--(\ref{saidameasured}) is dropped), the term $\Pi(k)$ can be removed.}
\end{remark}
       
\subsection{Global Convergence Result}

In the next theorem, we show that the proposed output-feedback controller based on monitoring function
drives $z$ to the $\Delta$-vicinity of the unknown maximizer $z^*$ defined in
\textbf{(A3)}. It does not imply that $z(t)$ remains in
$\mathcal{D}_\Delta$, $\forall t$. However, the oscillations
around $y^*$ can be made of order $\mathcal{O}(r)$.

\begin{theorem}                                                                           \label{theorem:01}
Consider the system (\ref{plant_inverse2})--(\ref{plant_extern2}), with output or cost function in (\ref{saidameasured}), control law (\ref{control_law}), modulation function (\ref{funcmod}), monitoring function (\ref{eqmonitpart})--(\ref{eqmonit}) and reference trajectory (\ref{modref}).
Assume that \textbf{(A1)--(A5)} hold, then: \textbf{(i)} the $\Delta$-vicinity  $\mathcal{D}_\Delta$ in \textbf{(A3)} is globally attractive being reached in finite time and \textbf{(ii)} for $L_{\Phi}$ sufficiently small, the oscillations around the maximum value $y^*$ of $y$ can be made of order $\mathcal{O}(r)$, with $r$ defined in (\ref{eqmonitpart}). Since the signal $y_m$ can be saturated in (\ref{modref}), all signals in a closed-loop system remain uniformly bounded.
\end{theorem}

\textbf{Proof:} Outside the $\Delta$-vicinity, the derivative of the cost function $\Phi(z)$ does not vanish ($d \Phi(z) /d z \neq 0\,, \forall z \notin \mathcal{D}_\Delta$). Thus, a lower norm bound $\underline{k}_p$ for $k_p$ can be obtained from the lower bound $L_{\Phi}$ given in  \textbf{(A3)} which is valid globally. Furthermore, Proposition \ref{prop:01} holds while $z$ stays outside the $\Delta$-vicinity, {\em i.e.}, no finite-time escape occurs for
the system signals.

 Now we proceed to the proof of the properties \textbf{(i)} and \textbf{(ii)} of Theorem~\ref{theorem:01}.

\textbf{(i)~Attractiveness of $\mathcal{D}_\Delta$}

This proof is made by contradiction. Assume that $z(t)$ stay outside the $\Delta$-vicinity for all $t$, i.e., $z \notin \mathcal{D}_\Delta$, $\forall t\in [0,t_M)$. Then, from Proposition $1$, there exists a finite time $t_s$ such that $\dot{e}=-\lambda_1\sgn(e)$, for some $\lambda_1>0$, $\forall t \geq t_s$. The error $e=y-y_m$ tends to zero, but since $y_m$ strictly increases with time and $y=\Phi(z)$ has a maximum value $y^*$, for $t$ large enough, $y_m>y^*\geq y$ and $\sgn(e)=-1$, assuring that $y$ increases with constant rate $(\dot y= k_m+\lambda_1)$, from (\ref{erro})-(\ref{modref}), that is, $y$ must approach $y^*$. So, $z$ is driven inside $\mathcal{D}_\Delta$, which is a contradiction. Thus, $\mathcal{D}_\Delta$ is attained in some finite time. Consequently, $z(t)$ remains or oscillates around $\mathcal{D}_\Delta$, and similarly $y$ with respect to some small vicinity of $y^*$, $\forall t$ large enough.

These oscillations come from the loss of control strength as $k_p \!\to\!0$ whenever the relation $\underline{k} _p\!\leq\!|k_p|$ is violated, or are due to the recurrent changes in the HFG sign at the extremum point $(z^*, y^*)$ where $k_p \!=\!0$ ($d \Phi(z)/dz=0$).  
In what follows, we show that these oscillations can be made ultimately of order $\mathcal{O}(r)$, with $r$ from (\ref{eqmonitpart}).
 
\textbf{(ii)~Oscillations of order $\mathcal{O}(r)$ around $y^*$}

According to the Assumption \textbf{(A3)}, $\Delta$ can be made arbitrary small so that $|y-y^*| = \mathcal{O}(r)$ when $z \in \mathcal{D}_{\Delta}$. Thus, if $z(t)$ remains in $\mathcal{D}_\Delta$, $\forall t$, the corresponding neighborhood of $y^*$ can be made of order $\mathcal{O}(r)$ with an appropriate $L_{\Phi}$. Otherwise, the key idea is to show that, if $z$ oscillates around $\mathcal{D}_\Delta$, this also holds since the time spent outside after leaving the set $\mathcal{D}_\Delta$ is also of order $\mathcal{O}(r)$.

Indeed, first reminding that after some finite time $t_{y^*}\!>\!0$, one has that $\sgn(e)\!=\!-1$, since $y_m$ strictly increases with time and $y=\Phi(z)$ has a maximum value $y^*$. Remind that the output error is given by
\begin{align}
	e(t)=y(t)-y_m(t), \quad \forall t>t_{y^*} \,.                                                                           \label{erro2}  
\end{align}

Now, assume that $z$ reaches the frontier of $\mathcal{D}_\Delta$ (from inside) at some time $T_{1}>t_{y^*}$ with wrong control direction at $t=T_1$. Note that $\mathcal{D}_\Delta$ is invariant when the control direction is correct and $t$ is large enough.

Moreover, from (\ref{erro2}), one can write
\begin{equation}
\tilde{e}(t)=\tilde{y}(t) -\delta_m(t-T_{1})\,,\label{errotil}
\end{equation}
where $\tilde{e}(t):=e(t)-e(T_{1})$, $\tilde{y}(t):= y(t)-y(T_{1})$, $\delta_m=0$ when $y_m$ is saturated and $\delta_m=k_m$, otherwise.

In addition, from (\ref{errotil}) one can also write
\begin{equation}                                                                                      \label{saidatil}    
\abs{\tilde{y}(t)} \leq \abs{\tilde{e}(t)}  + \delta_m(t-T_{1})\,.
\end{equation}

Let $T_2\!\geq\!T_1>t_{y^*}$ and $T_3\!\geq\!T_1>t_{y^*}$, where $T_{2}$ is the first time when $\varphi_m{(t)}=e(t)$ (independently if $z(t)$ is inside or outside $\mathcal{D}_\Delta$) and $T_3$ is the first time when $z(t)$ reaches the frontier of $\mathcal{D}_\Delta$ again (from outside). Notice that, $z(t)\!\notin\!\mathcal{D}_\Delta$ for $t\!\in\![T_{1},T_{3})$.

Now, consider two cases: {\bf (a)} $z$  reaches the frontier of $\mathcal{D}_\Delta$ with correct control direction ($T_3>T_2$) and {\bf (b)} $z$  reaches the
frontier of $\mathcal{D}_\Delta$ with wrong control direction ($T_3\leq T_2$).

For {\bf case (a)}, suppose $t\!\in\![T_1,T_3]$ and first consider $t\!\in\![T_1,T_2)$. During this time, the control law has wrong control direction. Thus, allowing $t_{y^*}$ to be large enough so that the exponential term $\pi_2/\underline{k}_{p}$  has been dominated by the term $\Pi(k)$ in the modulation function  $\rho$ in (\ref{funcmod}). 
Consequently, by the construction of the monitoring function, one has $e(T_1) < e(t) < e(T_2)$ and $e(T_{2})-e(T_1)<2r$. Otherwise, a change of control direction would occur 
according to Proposition \ref{prop:01}. Therefore, $\tilde{e}(t)=e(t)-e(T_{1})$ is of order $\mathcal{O}(r)$, $\forall t\!\in\![T_1,T_2)$. 
Moreover, since $\sgn(k_p)$ was wrongly estimated, one can write $|u(t)| = \rho(t)$, $\forall t\!\in\![T_1,T_2)$, with the control $u$ defined in (\ref{control_law}).

From (\ref{edoerro}), (\ref{moduerro}), (\ref{boundond}) and reminding that $t_{y^*}$ is considered large enough so that the exponential term $\pi_0$ has decreased to an arbitrary small value, one can verify that $|\dot{e}|\geq \underline{k}_p|u+d_e|$ and $|u+d_e|\geq \rho-|d_e|\geq\delta$ with an appropriate positive constant $\delta$ in (\ref{funcmod}).

Hence, $|\dot{e}(t)|>\delta_1$, $\forall t\!\in\![T_{1},T_{2})$, and $(t-T_1)\leq|\tilde{e}|/\delta_1$, where $\delta_1=\underline{k}_p\delta$ is an arbitrary positive constant.

Then, reminding that $\tilde{e}$ is of order $\mathcal{O}(r)$, $\forall t\!\in\![T_{1},T_{2})$, one can also assure that $(t-T_{1})$ and $\tilde{y}(t)$ in (\ref{saidatil}) are of order $\mathcal{O}(r)$, $\forall t\!\in\![T_{1},T_{2})$.

Moreover, by continuity, $\tilde{y}(t)$ is also of order $\mathcal{O}(r)$, $\forall t\!\in\![T_{1},T_{2}]$. Now, consider $t\!\in\![T_2,T_3]$. During
this interval, the control direction is correct$\sigma$ is in sliding motion, and thus $e\to 0$. From (\ref{erro2}) and (\ref{modref}), one has that $\dot{y} = k_m + \lambda_1 >0$ and $y(t)$ is strictly increasing, $\forall t \in [T_2,T_3]$. Consequently, one can conclude that $\tilde{y}(t)= y(t)-y(T_{1})$ is also of order $\mathcal{O}(r)$, $\forall t\!\in\![T_{2},T_{3}]$, since $y(t)$ approaches $y^*$ during the latter interval. Since this is also valid for the interval $[T_{1},T_{2})$, we have proved that the oscillation outside $\mathcal{D}_\Delta$ is of order $\mathcal{O}(r)$ in case (a), $\forall t\!\in\![T_{1},T_{3}]$.

For {\bf case (b)}, suppose $t\!\in\![T_1,T_3]$. During this time interval, the control law has wrong control direction. Thus, one can also conclude that $\tilde{y}(t)$ are of order $\mathcal{O}(r)$, $\forall t\!\in\![T_{1},T_{3}]$, following directly the steps of the first part of the proof of case (a).
By the continuity of the uncertain output function  $y=\Phi(z)$, the boundedness of $y$ stated above implies that $z$ is uniformly norm bounded
$(\mathcal{UB})$, and  also from \textbf{(A4)} one can easily conclude that all closed-loop system signals are  uniformly bounded. \hfill  $\qed$

\subsection{Multiple Extrema Envisage}
ESC applied to achieve global maximum in the presence of local extrema is a challenging area. 
Sometimes, the exhaustive search of the solution set may be the only choice, as discussed in \cite{TNMA:09}. The authors have presented a scalar extremum seeking feedback
controller that achieves semi-global practical global extremum
seeking despite local extrema.  

Inspired by the ideas introduced there, it was observed that by tuning a new design parameter in the monitoring function properly,  it is possible
to pass through a local extremum and converge to the global
one as well as is done in \cite{TNMA:09}, when the amplitude of the excitation (dither) signal is adaptively adjusted.

In this case,  the monitoring function  (\ref{eqmonitpart})--(\ref{eqmonit})  should be replaced by 
\begin{align} 
\varphi_k(t)= |e(t_k)|e^{-\lambda(t-t_k)}+ r + c(k)\,,                             \label{eqmonitpart2}
\end{align}   
where $c(k)$ is any positive monotonically decreasing sequence, such that $c(k)\to 0$ as $k \to +\infty$, and $c(0)>\mathcal{O}(r)$.

 \subsection{Illustrative example} \label{section4}

\textbf{Example 1} Consider in this example a plant with an unknown output performance characteristic and relative degree one dynamics described by 
\begin{align}                                                                   
	\dot{x}&=\left[
	\begin{array}{cc}
-1 & 1 \\ 1 &1		
	\end{array}\right]x+
		\left[
	\begin{array}{cc}
		0 \\ 1
	\end{array} \right]u	\\
		y=\Phi(z)&=e^{-\frac{(z-3)^2}{0.5}}+1.5e^{-\frac{(z-5)^2}{1.5}}\,, 
\end{align} 
where $x=[\eta^T \ z]^T$.

\begin{figure}[!htb]
\begin{center}
\includegraphics[width=.6\textwidth]{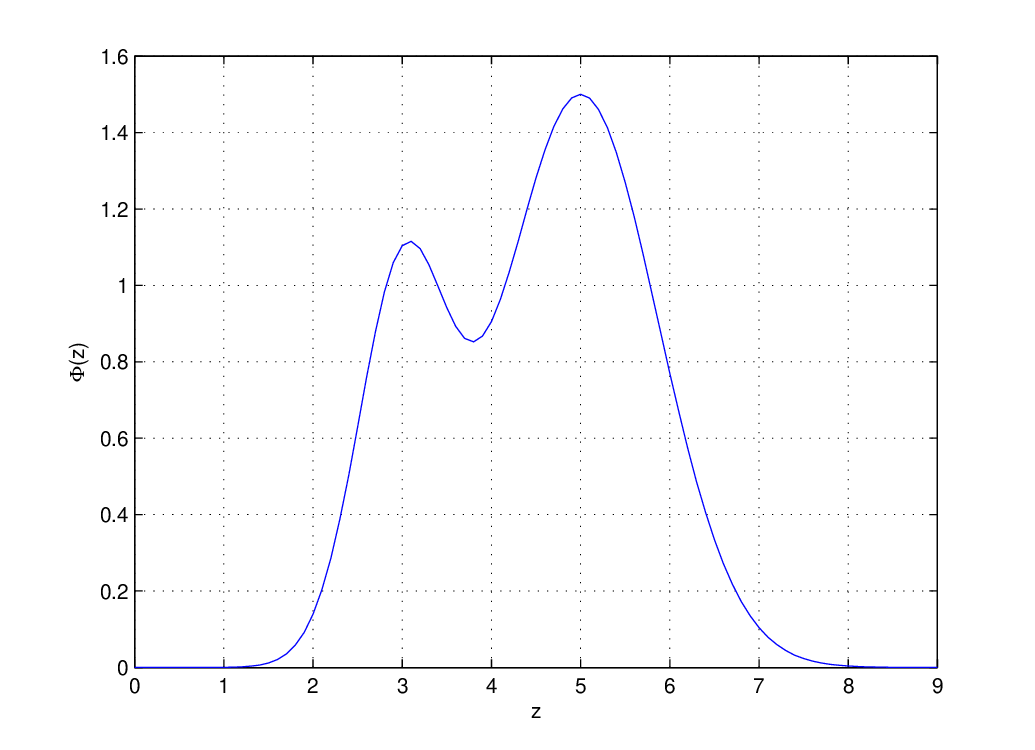}
\caption{Performance map $\Phi(z)$.}
\label{plantmultimax}
\end{center}
\end{figure}

The structure of this problem closely resembles various real-world goals, such as fine-tuning a one-dimensional engine calibration. In this scenario, the objective is to determine the optimal valve timing that maximizes efficiency \cite{NMMM:10}. Note that the performance map shown in Fig. \ref{plantmultimax} has multiple maxima, and, in contrast to \cite{NMMM:10}, global convergence (for all domain of initial conditions) can be guaranteed.

The modulation function was designed to satisfy (\ref{funcmod}) with $\Pi(k)=(k+1)e^{-t/(k+1)}$. The norm observer applied to determine $\bar\eta$ used in $\bar\phi_1(\bar\eta,z,t)$ has $\lambda_0=0.8$ and $\varphi_0(z,t)=2z$. The monitoring function (\ref{eqmonitpart2}) used to face the problem of local extrema has $c(k)=2/(k+1)$. 
The remaining parameters were: the lower bound $L_{\Phi}=\frac{2}{3}$, $\lambda=2$, $k_m=1$, $\delta=0.1$ and $r=0.1$.

\begin{figure}[!htb]
\begin{center}
\includegraphics[width=.55\textwidth]{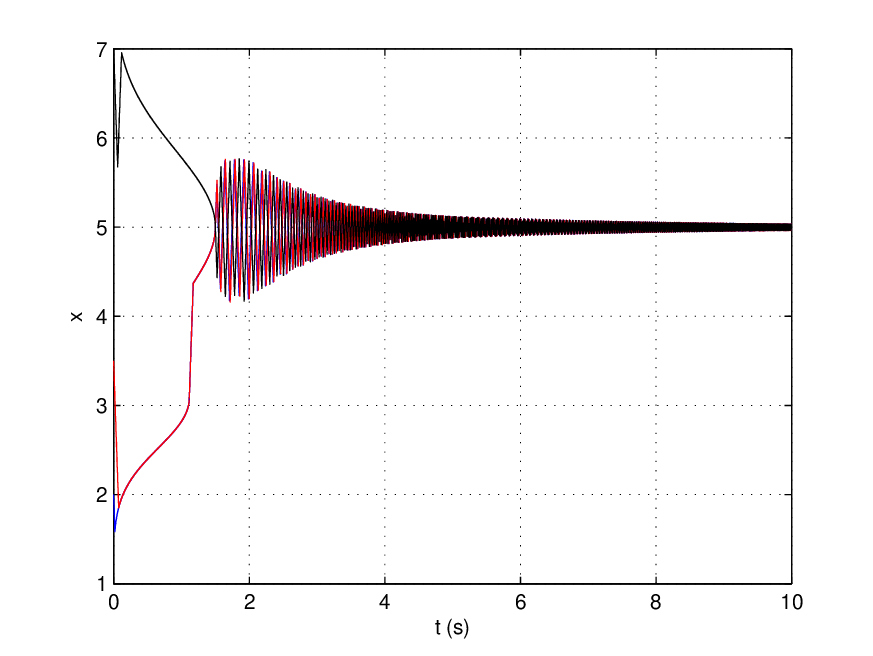}
\caption{Parameter $z$ converges to $z^*=5$  that maximizes $y$ using different initial conditions $z(0)$.}
\label{multimaxci2}
\end{center}
\end{figure}

\begin{figure}[!htb]
\begin{center}
\includegraphics[width=.55\textwidth]{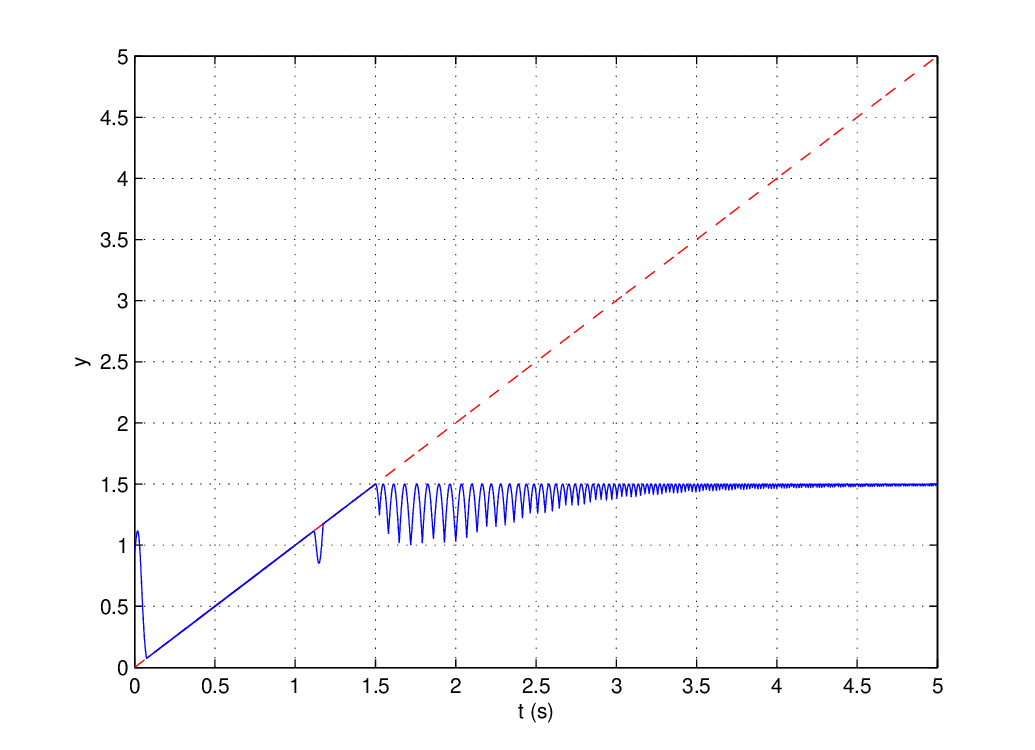}
\caption{Time history of the output plant $y$ (solid line) and the output model $y_m$ (dashed line). The output plant
tends to the maximum value $y^*=1.5$.}
\label{rastreamy_ym}
\end{center}
\end{figure}

\begin{figure}[!htb]
\begin{center}
\includegraphics[width=.55\textwidth]{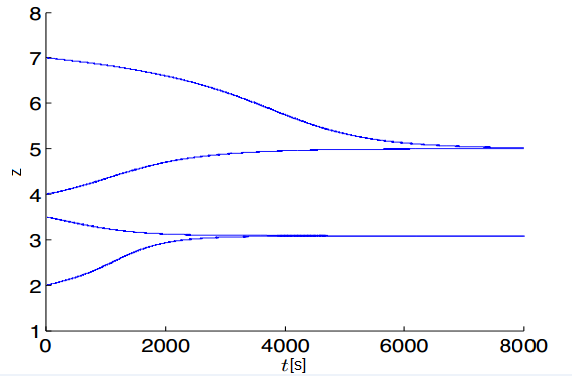}
\caption{Figure extracted from \cite{NMMM:10}. Convergence of $z$ to local maxima for different values of $z(0)$.}
\label{model_krstic}
\end{center}
\end{figure}

\begin{figure}[!htb]
\begin{center}
\includegraphics[width=.55\textwidth]{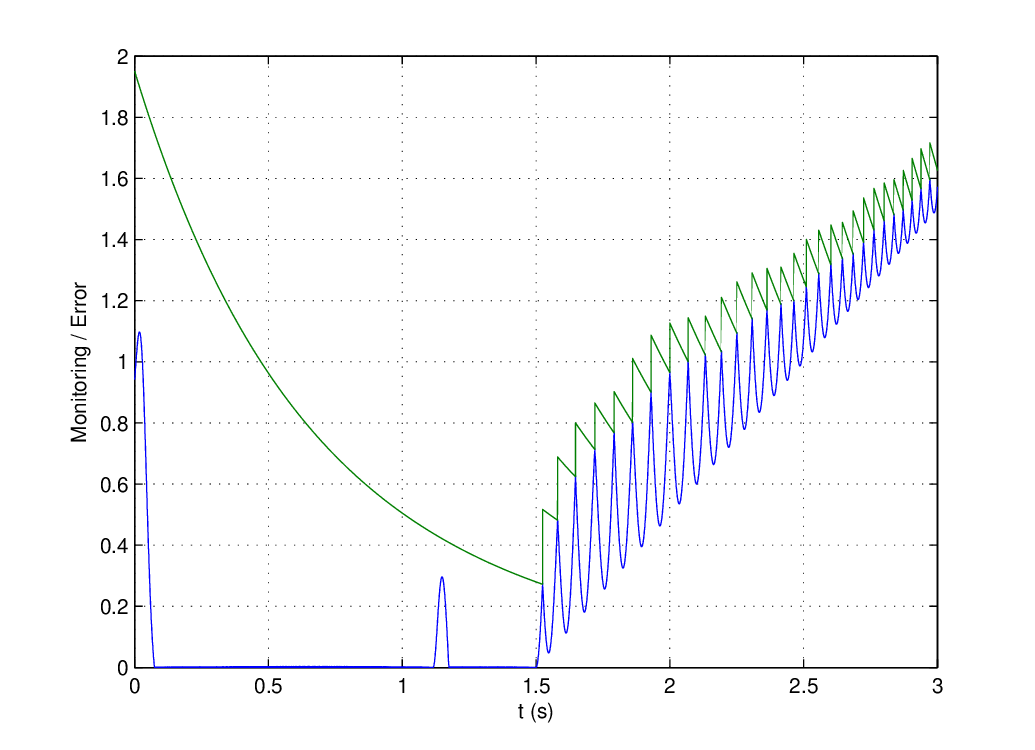}
\caption{Monitoring function $\varphi_m$ and error norm $|e|$.}
\label{funmonit_error}
\end{center}
\end{figure}

Fig. \ref{multimaxci2}  illustrates the convergence of the scheme for different initial conditions of $x$, corresponding to $z(0)=2$, $4$ and $7$. Note that, differently from \cite{NMMM:10} (see Fig. \ref{model_krstic}), where the convergence rate and global maximum are directly dependent on the initial conditions, here the example illustrates that it is possible to reach the global maximum in the presence of a local maximum without affecting the rate of convergence and independently of the initial conditions. 
As shown in Fig. \ref{rastreamy_ym}, $y$ tracks $y_m$ until $z$ reaches the vicinity of the global maximizer $z^*=5$. In Fig. \ref{plantmultimax}, the small oscillation just after $t=1$s shows the capacity of the monitoring function to pass through a local maximum at $y=1.1$ and converge to the global
one at $y=1.5$. Fig. \ref{funmonit_error} illustrates the monitoring function behavior. It can be seen that after reaching the global maximum, the error starts increasing because the reference trajectory is a ramp.

\break
\section{Uncertain and Arbitrary Relative Degrees} \label{uncertainrelativedegree}

Relative degree mitigation and the extremum searching are achieved
by means of a time-scaling technique. For the sake of simplicity,
we restrict ourselves to linear and stable plant dynamics
\cite{ZSK:2007}, but assuming that the relative degree can be
uncertain and arbitrary. By using a singular perturbation method,
it is shown that in the new time-scale, an attractive manifold is
revealed, which essentially reduces the considered system to a
single integrator perturbed by a fast sensor dynamics, which in
turn ultimately converges to a small residual set. We then exploit
this particular structure to redesign, with reduced control
authority, our original control law in Section~\ref{section3} to show its robustness with respect to the
arbitrary relative degree dynamics at the expense of time dilation,
which slows down the system response.

In this sense, we intend to show that the ESC proposed in
Section~\ref{section3} can also be extended to systems with uncertain and
arbitrary relative degree ($n^*$) in the form:
\begin{eqnarray}             
    \dot v &=& u \,,
    \label{eq0}\\
    \dot{x}&=& Ax+Bv \,,                                                                               \label{eq1new}  \\
    z&=&Cx \,,                         \label{outupunmeasurednew}
\end{eqnarray}
in cascade with a static subsystem
\begin{align}
y=\Phi(z)\,, \label{outputmensurednew}
\end{align}
where $u \in \mathbb{R}$ is the control input,  $x \in
\mathbb{R}^n$ is the state vector, $z \in \mathbb{R}$ is an
unmeasured output of the linear subsystem
(\ref{eq0})--(\ref{outupunmeasurednew}) and $y \in \mathbb{R}$ is
a measured output of the static subsystem
(\ref{outputmensurednew}), respectively.

\begin{remark} [Chattering Alleviation:]
\emph{The integrator in (\ref{eq0}) is used to obtain a virtual
control signal $v \in \mathbb{R}$, which increases the relative
degree of the system \cite{Lev:03}, \emph{i.e.}, $n\geq n^*-1$
instead of $n\geq n^*$. The increase of the relative degree yields
high-frequency switching in the control $u$, while the virtual
control $v$ driving directly the plant is continuous. Thus,
chattering attenuation is expected to be improved \cite{UGS:99}.}
\end{remark}

The matrices $A\in\mathbb{R}^{n\times n}$, $B\in\mathbb{R}^{n}$,
$C\in\mathbb{R}^{1\times n}$ and the order $n$ of the subsystem
(\ref{eq1new}) may also be uncertain.
The uncertain nonlinear function $\Phi:\mathbb{R}\to\mathbb{R}$ to
be maximized must still satisfy the Assumption $\bf{(A1)}$. The matrix $A$ in
(\ref{eq1new}) must be Hurwitz, to guaranty the uncertainty bounds
for control design.

\subsection{Singular Perturbation Analysis}
\label{coco}

We start summarizing some results obtained in Section~\ref{section2} for systems with
relative degree one. Consider the simplest case of the
integrator below with a nonlinear output mapping:
\begin{eqnarray}
  \dot z &=& u\,,                                              \label{unitary_plant_inverse2}\\
  y&=&\Phi(z)\,,                                            \label{unitary_saidameasured}
  \end{eqnarray}
where $u \in \mathbb{R}$ is the control input,  $z \in \mathbb{R}$
is the state vector and $y \in \mathbb{R}$ is the measured output
of the static subsystem.
In order to present such generalization, consider the system
(\ref{unitary_plant_inverse2})--(\ref{unitary_saidameasured}),
with the change of variables $z=v$, such that
\begin{eqnarray}
  \dot v &=& u\,,                                              \label{unitary_plant_inverse3}\\
  y&=&\Phi(v)\,,                                               \label{unitary_saidameasured3}
  \end{eqnarray}
can be directly controlled by the method of monitoring function
described in Section~\ref{monitsection}.

By using the singular perturbation approach \cite{KKO:1986}, it
can be shown that the monitoring function method for ESC
(see Section~\ref{section2}) is robust to fast
unmodeled dynamics such that the perturbed system
(\ref{unitary_plant_inverse3})--(\ref{unitary_saidameasured3}) is
rewritten in the following \emph{block sensor form} \cite[p.
50]{KKO:1986}
\begin{eqnarray}             
    \dot v &=& u \,,                               \label{fast_eq0}\\
    \mu\dot{x}&=& Ax+Bv \,,                                                                               \label{fast_eq1}  \\
    y &=& \Phi(C x) \,,                                                                                                \label{fast_outupunmeasured}
\end{eqnarray}
and ultimately satisfies the inequality
\begin{eqnarray}             
    |y-y^*| \leq \mathcal{O}(\sqrt{\mu})\,,    \label{residual_set}
\end{eqnarray}
where 
$\mu>0$ is a sufficiently small constant. The
complete demonstration of (\ref{residual_set}) follows the same
steps
presented in \cite{CH:1991,CH:1992}. 

In the singular case $\mu=0$, the differential equation
(\ref{fast_eq1}) is replaced by the algebraic equation
$x=-A^{-1}Bv$ and, from (\ref{fast_eq0}) and
(\ref{fast_outupunmeasured}), the first time derivative of the
output signal $y$ is given by
\begin{align}
\dot{y} = k_p(z) u\,,\label{eq:ydynamics_singular}
\end{align}
where the HFG is now rewritten as
\begin{align}
k_p(z)=-\Phi'(z)CA^{-1}B\,. \label{kpneras_singular}
\end{align}

From (\ref{kpneras_singular}) and $(\bf{A3})$, $k_p(z)$ satisfies
%
\begin{equation}
|k_p(z)|\geq \underline{k}_{p}>0\,, \quad \forall z \notin
\mathcal{D}_{\Delta}\,, \label{kpbarneras_singular}
\end{equation}
where $\underline{k}_{p}\leq L_\Phi |CA^{-1}B|$ is a known
constant lower bound for the HFG, considering all the admissible
uncertainties in $\Phi(\cdot)$, $A$, $B$, and $C$.

\subsection{Time-Scaling for Control Redesign}

Thus, by applying an appropriate linear time-scaling
\cite{MONPP:2002}
\begin{equation}
  \frac{d t}{d \tau} = \mu\,,                                              \label{time-scaling}\\
\end{equation}
the system (\ref{fast_eq0})--(\ref{fast_outupunmeasured})  can be rewritten as 
\begin{eqnarray}
  v' &=& \mu u                               \label{eq02}\\
  x' &=& A x + Bv\,,                                              \label{plant_inverse2new}\\
  z &=& C x \,,                     \label{plant_extern2new}\\
  y&=&\Phi(z)\,,                                                                          \label{saidameasurednew}
  \end{eqnarray}

where $v':=\frac{d v}{d \tau}$ and $x':=\frac{d x}{d \tau}$.
It means that $\exists \mu^*>0$ such that the input signal $u$ can
be scaled (\ref{eq02}) to control the original system
(\ref{eq1new})--(\ref{outputmensurednew})
in a different dilated time-scale governed by $t=\mu\tau$,
$\forall \mu \in (0,\mu^*]$.

The significance lies in the fact that the monitoring function-based ESC, initially designed for systems with a relative degree of one, demonstrates robust to rapid unmodeled stable dynamics ($\mu\to +0$). Consequently, it proves suitable for controlling dynamics of arbitrary relative degree when appropriately scaled. However, it's anticipated that as $\mu\to +0$, the trade-off is a slower response of the closed-loop system.

\subsection{Scaled Controller Parameters}

The time-scaling (\ref{time-scaling}), allows us to consider the
original plant (\ref{eq1new})--(\ref{outputmensurednew}) in a
different time-scale being controlled by the controller
(\ref{control_law}) properly scaled by $\mu u$, see (\ref{eq02}).

In order to incorporate it, the modulation function must be
redesigned to satisfy
\begin{align}
\rho\geq\mu[|d_e|+\delta], \label{fast_funcmodgen}
\end{align}
instead of (\ref{funcmodgen}).

Remind that the derivative of the objective function does not
vanish $\forall z$ outside the $\Delta$-vicinity. Thus, the lower
norm bound $\underline{k}_p$ for $k_p$
given in (\ref{kpbarneras_singular}) holds.

Therefore, one can obtain the following norm bound as $d_e$
defined in (\ref{moduerro}):
\begin{equation}                                                                                          \label{boundond2}
|d_e(t)|\leq \bar{d}_e\,, \quad \bar{d}_e:=(k_m+\lambda
|e|)/\underline{k}_{p}\,.
\end{equation}

In the proposed scheme,  the following proposition provides one
possible modulation
function implementation so that (\ref{fast_funcmodgen}) holds. 

\begin{proposition}                                                                      \label{prop:01cdc14}
Consider the system (\ref{eq0})--(\ref{outputmensurednew}), the
reference model (\ref{modref}) and control law
(\ref{control_law}). Outside the $\Delta$-vicinity
$\mathcal{D}_\Delta$, if $\rho$ in (\ref{control_law}) is designed
as
\begin{equation}
\rho:=\frac{\mu}{\underline{k}_p} \left[k_m+\lambda |e|\right] +
\mu\delta\,, \label{funcmod2}
\end{equation}
for $\mu$ sufficiently small, then, while $z \notin
\mathcal{D}_\Delta$, one has: \textbf{(a)} the switching of the monitoring function (\ref{eqmonitpart})--(\ref{eqmonit}) stops, \textbf{(b)} no finite-time escape of the system
signals occurs, which implies that $t_M \to +\infty$,  and
\textbf{(c)} the error $e(t)$ tends to zero in finite time.  The
term $\delta>0$ is any arbitrary constant.
\end{proposition}
\textbf{Proof:} By considering the singular perturbation argument
and the time-scaling (\ref{time-scaling}), which show that the
systems (\ref{fast_eq0})--(\ref{fast_outupunmeasured}) and
(\ref{eq02})--(\ref{saidameasurednew}) are equivalent for $\mu$
sufficiently small, then, the demonstration for the original plant
(\ref{eq0})--(\ref{outputmensurednew}) follows the same steps
presented in the proof of
Proposition ~\ref{prop:01} for the relative degree one case.  \hfill  $\qed$

\begin{remark} [Tuning Rules:]
\emph{Since the control design is developed in the light of the
slow time-scale $\mu t$, it is natural that the parameters $k_m$
of the reference model (\ref{modref}) and $\lambda$ of the
monitoring function (\ref{eqmonitpart})--(\ref{eqmonit}) are
appropriately scaled, \emph{i.e.}, replaced by $\mu k_m$ and $\mu
\lambda$, respectively. In our ESC application, the ultimate
residual set of the proposed algorithm around the maximum $y^*$ is
dependent on the values to which the monitoring function
converges. According to the definition given in
(\ref{eqmonitpart}), the ultimate residual set $r$ must be of
order $\mathcal{O}(\sqrt{\mu})$ to satisfy (\ref{residual_set}).}
\end{remark}

\subsection{Global Convergence Result}

The following theorem states that the proposed output-feedback
controller based on monitoring function drives $z$ to the
$\Delta$-vicinity defined in \textbf{(A1)} of the unknown
maximizer $z^*$. It does not imply that $z(t)$ remains in
$\mathcal{D}_\Delta$, $\forall t$. However, the amplitude of
signal oscillations around $y^*$ can be kept of order
$\mathcal{O}(\sqrt{\mu})$.


\begin{theorem}                                                                           \label{theorem:01cdc14}
Consider the system (\ref{eq1new})--(\ref{outupunmeasurednew}),
with output or objective function in (\ref{outputmensurednew}),
control law (\ref{eq0}) and (\ref{control_law}), modulation
function (\ref{funcmod2}), monitoring function
(\ref{eqmonitpart})--(\ref{eqmonit}) and reference model (\ref{modref}). 
Assume that \textbf{(A1)--(A3)} hold, then: \textbf{(i)} the
$\Delta$-vicinity $\mathcal{D}_\Delta$ in \textbf{(A1)} is
globally attractive being reached in finite time and \textbf{(ii)}
for $L_{\Phi}(\mu)$ in \textbf{(A1)} sufficiently small, the
oscillations of $y(t)$ around the maximum value $y^*$ can be made
of order $\mathcal{O}(\sqrt{\mu})$, with $r$ defined in
(\ref{eqmonitpart}) being also a constant of order
$\mathcal{O}(\sqrt{\mu})$. Since the signal $y_m$ can be saturated
in (\ref{modref}), all signals in closed-loop system remain
uniformly bounded.
\end{theorem}

\textbf{Proof:} As in the proof of Proposition~\ref{prop:01cdc14}, the
demonstration is based on singular perturbation/time-scaling
arguments and follows the steps presented in the proof of
Theorem~\ref{theorem:01}, for $\mu$ sufficiently small. \hfill  $\qed$

\section{Light Source Seeking Experiments} \label{section4cdc14}

\begin{figure}[!htb]
\begin{center}
\includegraphics[width=.46\textwidth]{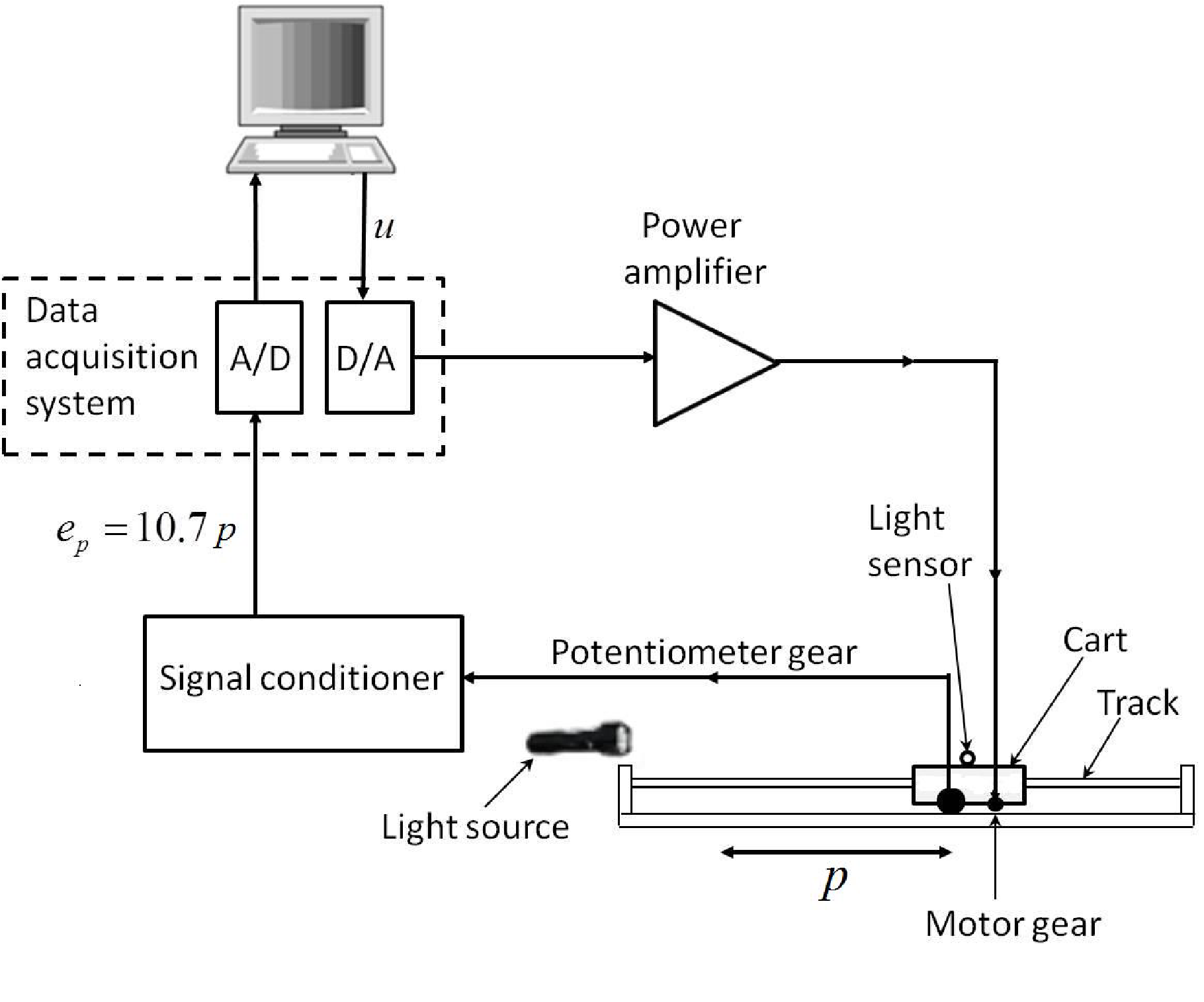}
\caption{Block diagram of the source seeking experiment.}
\label{plantmultimax2}
\end{center}
\end{figure}

\begin{figure}[!htb]
\begin{center}
\includegraphics[width=.6\textwidth]{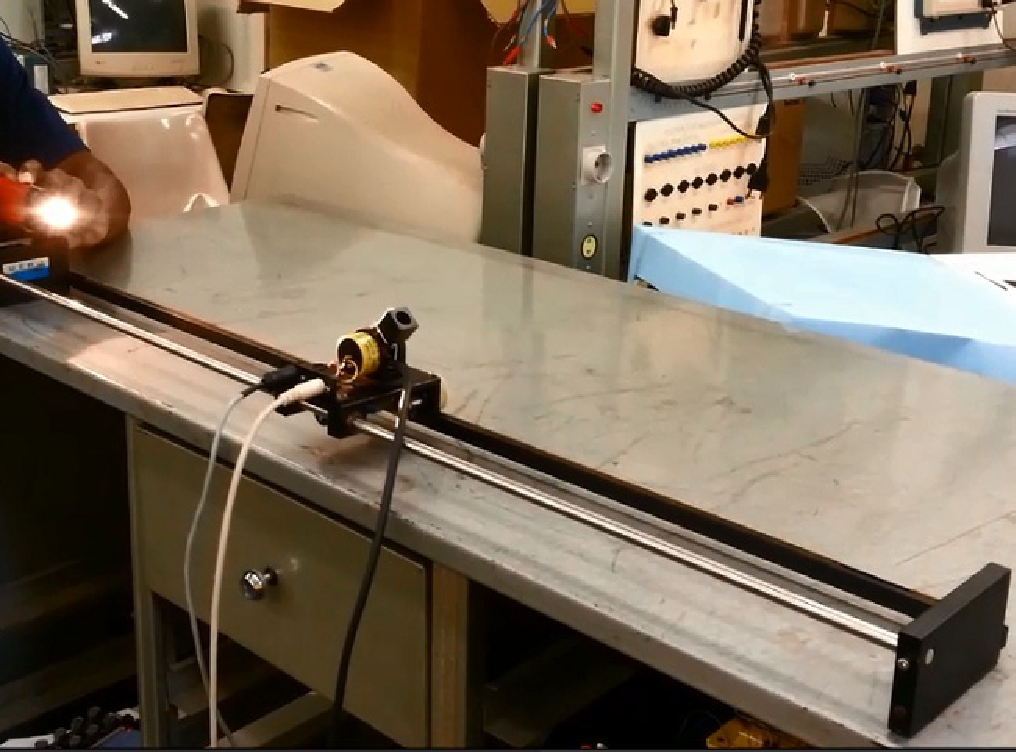}
\caption{Experiment~$I$: fixed light source.} \label{setup1}
\end{center}
\end{figure}

\begin{figure}[!htb]
\begin{center}
\includegraphics[width=.6\textwidth]{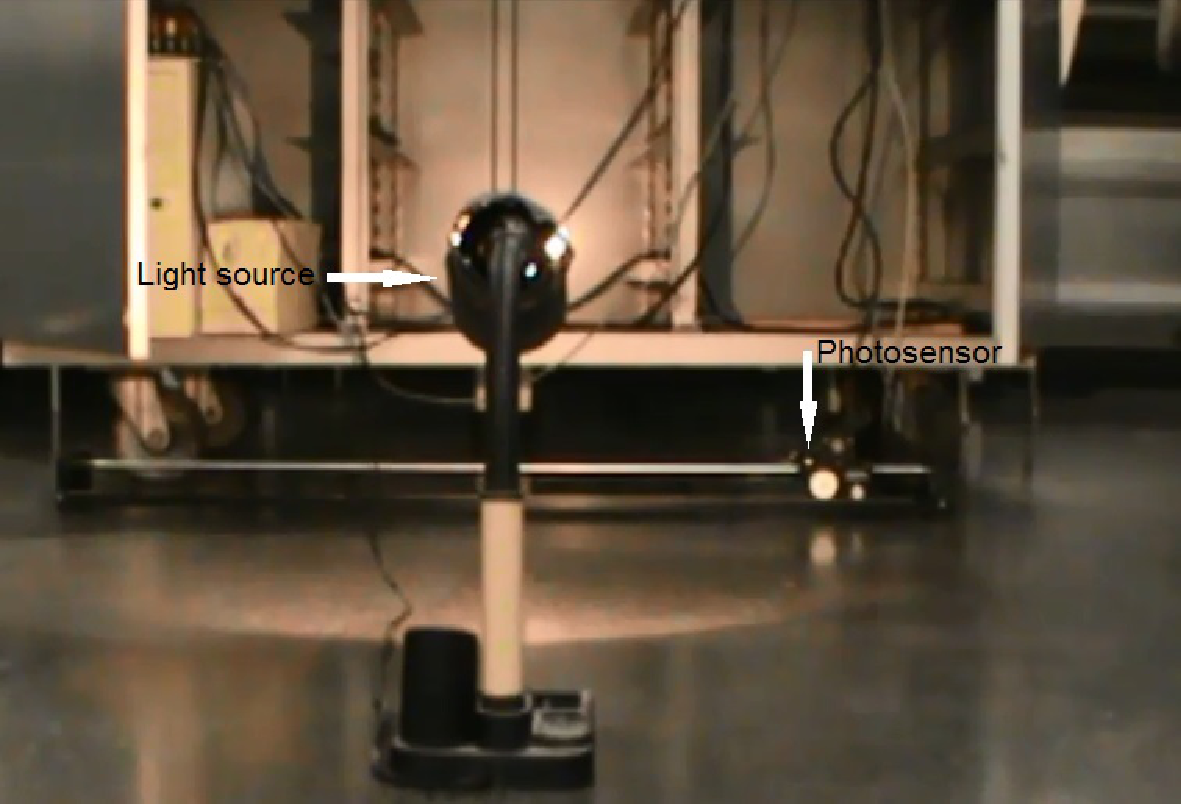}
\caption{Experiment~$II$: moving light source.} \label{setup2}
\end{center}
\end{figure}

The basic components utilized in the experiments involving linear motion comprise the cart and track are depicted in Figs. \ref{plantmultimax2}, \ref{setup1}
and \ref{setup2}, provided by Quanser Consulting (Linear Position Servo module IP01). The setup includes a cart that smoothly glides along a stainless steel shaft on the ground. Within the cart, a DC motor and a potentiometer are integrated. These components are interconnected through a rack and pinion mechanism: the DC motor contributes the driving force to the system, while the potentiometer facilitates the measurement of the cart's position.
The motor shaft is connected to a gear with a diameter of $0.5"$, and the potentiometer shaft is connected to a gear with a diameter of $1.17"$. Both gears engage with a toothed rack. As the motor rotates, the torque generated at the output shaft is transformed into a linear force, initiating the movement of the cart. Concurrently, as the cart advances, the potentiometer shaft undergoes rotation, and the voltage, denoted as ($e_p$), measured from the potentiometer can be calibrated to determine the position ($p$), \emph{i.e.}, $e_p = 10.7p$, along the track.

The control algorithm was implemented in a digital computer using
the Euler integration method with $1ms$ step size. The car
position ($p$) is measured by a potentiometer connected to a $12$
bit {A/D} converter through a signal conditioner, which reduces
aliasing and keeps measurement noise amplitude $<0.4 mm$
peak-to-peak. The control signal ($u$) is a voltage generated by a
$12$ bit {D/A} converter connected to a linear power operational
amplifier, which drives the {DC} motor. The saturation limits the
control signal peak amplitude to $10V$.

The objective of our feedback system for source-seeking is to regulate the cart's position without relying on data from the potentiometer. Instead, it utilizes the light intensity received by a photosensor, which is also mounted on the top of the cart, as illustrated in Fig.~\ref{setup2}. The output from the photosensor undergoes filtration through a low-pass RC filter to eliminate high-frequency noise. It's important to note that the position signal derived from the potentiometer will solely serve the purpose of monitoring system variables in our experiments, without contributing to the design of the control. Simply put, the control designer lacks information about both the cart's position and the source position.

\begin{remark}[Uncertain Relative Degree:]
\emph{ The breakdown of a controlled system into distinct components: actuator, plant, and sensor is unique. Take, for instance, the flexibility to integrate any actuator or sensor directly into the plant. This integration can significantly alter the relative degree \cite{L:2012}, underscoring the importance of the chosen control strategy. A notable drawback of conventional sliding mode-based control methods is their reliance on a precisely defined, constant, and known relative degree for the sliding variable. Even minor disturbances or inaccuracies in the model can cause a reduction in the relative degree or its complete elimination
\cite{L:2012,H:2001}.}
\end{remark}

In what follows, we evaluate the effectiveness of the proposed
controller based on monitoring function and time-scaling under
perturbations of external lighting conditions (the light of the
laboratory was not turned off during the essays). The performance
has also been shown to be robust with respect to the presence of
the actuator and sensor dynamics \cite{GK:2011}, which can be
included in the design stage as a part of the linear subsystem
(\ref{eq1new}).

The {DC} motor of our experimental scenario is essentially
modelled as a linear plant with the following nominal first order
transfer function:
\begin{equation}\label{iosystem}
G_p(s)=\frac{z}{v}=\frac{3.9}{(s+17.2)}\,,
\end{equation}
where $z$ is the angular velocity in $rpm$, the armature voltage
$v$ is the control input and the model parameters are considered
uncertain in the control design.

The state-space equation of the overall system considering the
integrator (\ref{eq0}) for chattering alleviation, the relative
degree one dynamics (\ref{iosystem}) of the {DC} motor and the
unknown nonlinear field $\Phi(\cdot)$ created by the light source
can be written as 
\begin{align}\label{realsystem}
    \dot{v}&=u\\
    \dot{z}&=-17.2 z + 3.9 v \\
        y&=\Phi(z)\,,\label{merdaaoquadrado}
\end{align}
where only the output $y$ measured by the photosensor is
assumed available for feedback in our experiments. The maximum of
$\Phi(z)$ occurs near of the light source, \emph{i.e.}, where the
gradient of the field is practically null. The main idea for the
control law is to guide the cart up the gradient of the signal to
find the source.

Two experiments were conducted in this study. The primary experiment aimed to confirm the effectiveness of the proposed controller in localizing a stationary light source (see
Fig.~\ref{setup1}). The secondary experiment focused on assessing the tracking capabilities of the proposed ESC when the light source was in motion (Fig.~\ref{setup2}).

\begin{figure}[!htb]
\begin{center}
\includegraphics[width=.6\textwidth]{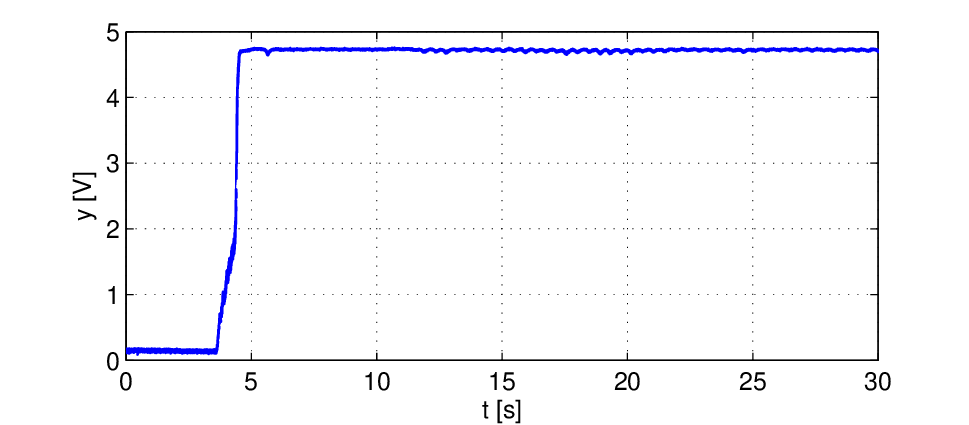}
\caption{Moving source: light intensity received by the
photosensor.} \label{multimaxci3}
\end{center}
\end{figure}

In both experiments, the reference signal was chosen as in
(\ref{modref}) with $k_m=2\mu$ and $y_m$ was saturated at $5$
since this was the maximum light intensity given by the circuit of
the photosensor, see Fig.~\ref{multimaxci3}. Consequently, $e$ is
uniformly bounded.

The monitoring function (\ref{eqmonitpart}) used to solve our
source seeking problem has $\lambda=\mu$ and $r=0.2{\sqrt\mu}$.
The modulation function in the control law (\ref{control_law}) was
designed to satisfy (\ref{funcmod2}).      In this case, the controller parameters could be:
$\underline{k}_p=0.2 L_{\Phi}$, $L_{\Phi}=20 r$, $\delta=0.1$ and
$\mu=0.5$.                                                                       \label{funmod2}

\begin{figure}[!htb]
\begin{center}
\includegraphics[width=.6\textwidth]{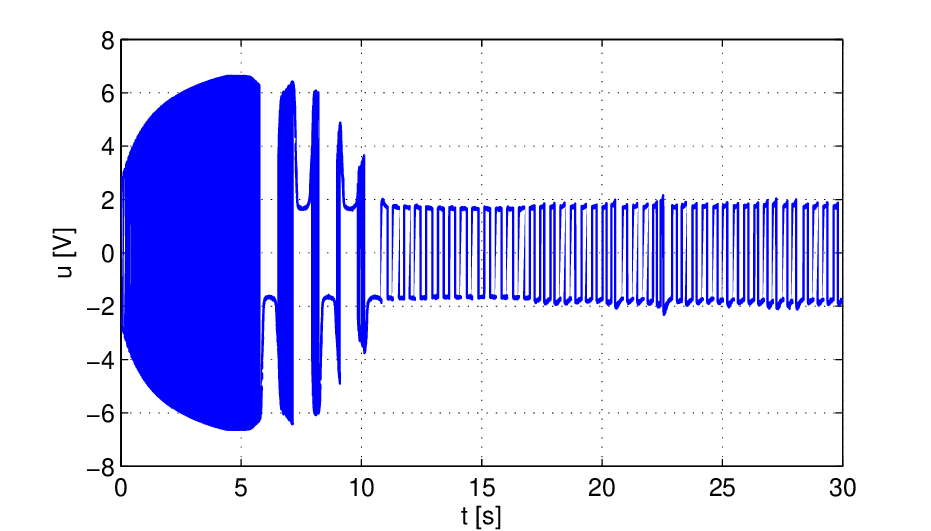}
\caption{Control signal $u(t)$.} \label{rastreamy_ym2}
\end{center}
\end{figure}

The discrete-time implementation of the controllers and unmodeled
dynamics induce the chattering phenomenon \cite{UGS:99} seen in
the control signal (Fig.~\ref{rastreamy_ym2}). 
However, it is not a problem since the control signal $v$ applied
to the plant is smooth (curve not shown) due to the integral
action in (\ref{realsystem}) used to reduce control
chattering, see Remark~2. 

In what follows, we discuss in detail only the results of
Experiment II. In this experiment, the light source is positioned
in front of the track at a distance of approximately $1m$
(Fig.~\ref{setup2}). The ESC moves the cart along of the track to
search the point of maximum light power.

\begin{figure}[!htb]
\begin{center}
\includegraphics[width=.6\textwidth]{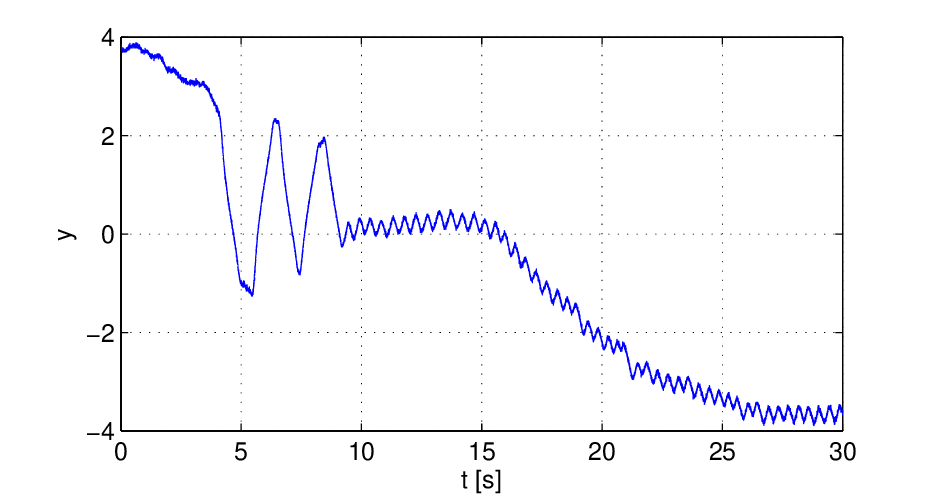}
\caption{Cart position $p(t)$ along of the track.}
\label{funmonit_error2}
\end{center}
\end{figure}

Fig.~\ref{multimaxci3} presents the corresponding output voltage
of the light intensity $y$ received by the photosensor. Notice
that in the initial response up to $t=4s$, the light source was
turned off. The residual set captured within this time interval is
due to the ambient light. After that, the light source is turned
on and the measurement keeps almost constant even when the source
is moved and the cart tries to maximize the light intensity $\Phi$
which is being measured, see Fig.~\ref{funmonit_error2}.

Fig.~\ref{rastreamy_ym2} illustrates the temporal evolution of the control signal. Prior to the activation of the light source, there is noticeable high-frequency switching of the control direction due to a violation of an inequality condition (\ref{eq:boundMIMOKpconhecido2}). This phenomenon arises because the output error $e(t)$ tends to increase while the reference signal $y_m(t)$ saturates at 5 volts, leading to a larger error signal during the initial phase, as indicated in Fig.~\ref{model_krstic}. Additionally, Fig.~\ref{model_krstic} reveals a transition to lower frequency switching of the control direction between t = 0 and 10 seconds, coinciding precisely with the detection of the maximum power point of the light source. Subsequently, the cart settles into oscillations around this point. Furthermore, Fig.~\ref{model_krstic} provides insight into the behavior of the monitoring function $\varphi_m(t)$ and the output error norm during the movement of the light source, as described earlier.

\begin{figure}[!htb]
\begin{center}
\includegraphics[width=.6\textwidth]{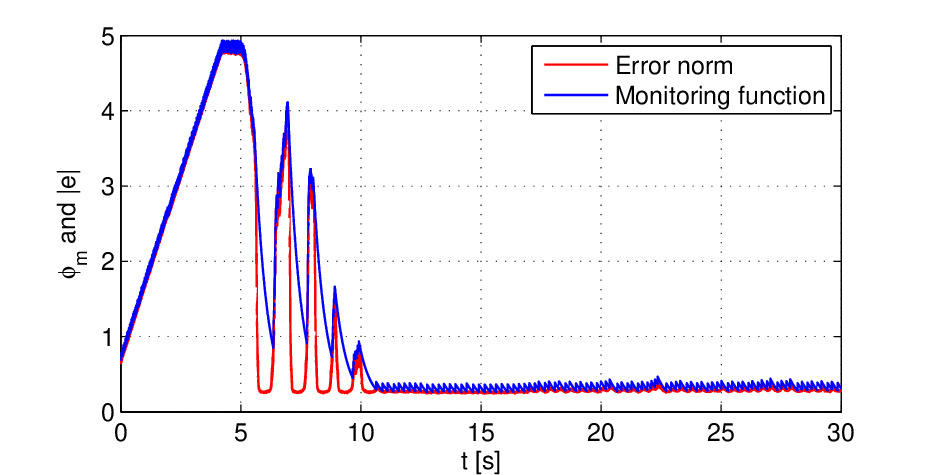}
\caption{Time history of $\varphi_m(t)$ and $|e(t)|$ when the
light sensor is moving during the experiment.}
\label{model_krstic2}
\end{center}
\end{figure}

Fig.~\ref{funmonit_error2} illustrates the cart position during the
experiments. For $t\in[0,4]$, the cart is practically stopped.
Immediately after the photosensor has detected the presence of the
light source, the cart is moved into its direction and oscillates
with decreasing amplitudes around the light source until to reach
the ultimate residual set of order $\mathcal{O}(r)$. Since the ESC
algorithm never stops searching, the cart continues to sniff
around the light source. For $t\in[15,30]$, the light source is
slowly moved in parallel along of the track and is quickly
followed by the cart. Surprisingly, the light intensity and the
amplitude of the oscillations are kept practically constant as
shown in Figs.~\ref{multimaxci3} and \ref{model_krstic2},
respectively.

Fig.~\ref{funmonit_error2} depicts the position of the cart throughout the experiments. From $t=0$ to $4$ seconds, the cart remains nearly stationary. Once the photosensor detects the light source, the cart promptly moves towards it, initiating oscillations with diminishing amplitudes around the light source until reaching a final residual set point denoted as $\mathcal{O}(r)$. As the ESC algorithm continuously searches, the cart persists in probing around the light source. Between $t=15$ and $30$ seconds, the light source gradually shifts parallel along the track, swiftly trailed by the cart. Remarkably, both the light intensity and oscillation amplitude remain practically constant, as evidenced in Figs.~\ref{multimaxci3} and \ref{model_krstic2}, respectively.

The straightforward servomechanism facilitated a concise delineation of the controller design and enabled the assessment of closed-loop performance in real experimental scenarios involving disturbances (such as dry friction, varying lighting conditions), measurement noise, unmodeled dynamics, and significant uncertainties.

The videos of the experiments can be found in the following link:
{\color{blue}{\url{https://bitily.me/eTUXE}}}.

\section{Conclusions}
A novel extremum seeking control scheme, incorporating monitoring
functions, norm state estimation and time-scaling has been derived for a class of
uncertain nonlinear plants. This approach ensures the global convergence
of the system output to a small neighborhood of the extremum
point using only output feedback. Simulation results have been
conducted to illustrate the controller performance,
even in scenarios involving local extrema. To our knowledge,
real-time solutions with global convergence properties, exclusively
based on output feedback, did not exist. This paper addresses this
gap and introduces a solution. Furthermore, we establish that our
extremum seeking controller, based on monitoring functions, touches on engineering applications, including the localization
and tracking of a light source. The theoretical results are evaluated  through successful experimental tests in a one-dimensional source seeking
problem. This problem involves the control of a linear
servomechanism cart without the use of position/velocity measurements. The generalization of the proposed method to higher-order real-time optimization schemes involving gradient and Hessian estimates \citep{REOMK:19,OK:2015,GKN:12b} of the map rather than using only its output measurement seems to be an interesting topic for future investigation.

\section*{Aknowlegments}

This study was financed in part by the Coordenação de Aperfeiçoamento
de Pessoal de Nível Superior - Brasil (CAPES) Finance
Code 001; Conselho Nacional de Desenvolvimento Científico e Tecnológico,
CNPq; Fundação de Amparo à Pesquisa do Estado do Rio
de Janeiro, FAPERJ.

\end{document}